\documentstyle[12pt]{article}
\begin{document}

\title{On geometry and mechanics\footnote{Dedicated to the memory of my grandparents Nikolaos and Alexandra, and Konstantinos and Eleni.}}

\author{Nikolaos E. Sofronidis\footnote{$A \Sigma MA:$ 130/2543/94}}

\date{\footnotesize Department of Economics, University of Ioannina, Ioannina 45110, Greece.
(nsofron@otenet.gr, nsofron@cc.uoi.gr)}

\maketitle

\begin{abstract}
Our purpose in this article is first, following [14], to find the topological upper limits of projections of secant planes to $C^{1}$ surfaces and the topological upper limits of projections of secant hyperplanes to $C^{1}$ hypersurfaces and second to prove that every $C^{2}$ space
curve can be the solution of the principle of stationary action.
\end{abstract}

\section*{\footnotesize{{\bf Mathematics Subject Classification:} 03E15, 03E65, 26B99, 49J10, 53A05, 70H03.}}

The definition of a manifold tangent vector for dimensions $5$,
$6$, $7$, ... (see [17]) makes use of the Axiom of Choice and hence
is a theorem of at least {\it ZFC - Axiom of Foundation}, which
implies the Banach - Tarski paradox (see [18]), since for any point
$P$ and for any tangent vector $v$ at $P$ one has to choose a path
$c$ such that $v = c'(0)$, while the Banach - Tarski paradox
contradicts the Ideal Gas Law (see [4] or [19]). Since in dimensions $2$,
$3$, $4$ the projection of the intersection of the secant
hyperplane with the hypersurface is in dimensions $1$, $2$, $3$,
one has the classical notion of tangent vector. So our purpose in
this article is to prove in {\it ZF - Axiom of
Foundation + Axiom of Countable Choice} two theorems for the property of being
tangent and an exclusion principle for the potential in Lagrangian mechanics.

\section{Two theorems on geometry}

{\bf 1.1. Definition.} If $X$ is any compact Polish space, then we
denote by $K(X)$ the compact Polish space of compact subsets of
$X$ equipped with the Hausdorff metric and for any sequence
$\left( K_{n} \right) _{n\in {\bf N}}$ of compact subsets of $X$,
we denote by $\limsup\limits_{n \rightarrow \infty }K_{n}$ the
topological upper limit of $\left( K_{n} \right) _{n\in {\bf
N}}$, i.e., the set of all $x \in X$ with the property that for
any $n \in {\bf N}$, there exists $x_{n} \in K_{n}$ such that $x$
is a limit point of the sequence $\left( x_{n} \right) _{n\in
{\bf N}}$. (See, for example, Section 4.F on pages 24-28 of [12]
and page 18 of [2] as it is cited on page 341 of [10].)
\\ \\
{\bf 1.2. Definition.} Let $\Omega \in K \left( {\bf R}^{2} \right)
$ be such that ${\Omega }^{ \circ } \neq \emptyset $ and let $f :
\Omega \rightarrow {\bf R}$ be $C^{1}$, while $\left( x_{0}, y_{0}
\right) \in {\Omega }^{ \circ }$. If $A$, $B$ are any real numbers
and $$P(A,B) = \left\{ (x,y,z) \in {\bf R}^{3} : z - f \left(
x_{0}, y_{0} \right) = A \left( x - x_{0} \right) + B \left( y -
y_{0} \right) \right\} $$ is any plane in ${\bf R}^{3}$ that
passes through $\left( x_{0}, y_{0}, f \left( x_{0}, y_{0} \right)
\right) $, then $P(A,B)$ intersects the surface $Graph(f)$ whose
equation is $z = f(x,y)$ exactly at the points $(x,y,z)$ in ${\bf
R}^{3}$ that satisfy the equation $f(x,y) - f \left( x_{0}, y_{0}
\right) = A \left( x - x_{0} \right) + B \left( y - y_{0} \right)
$ (see, for example, [9] and Paragraphs a and b of Section 8.4.6 on pages
470-471 of [16] and pages 495-497 of [11] and pages 326-327 of [3])
and let $$C(A,B) = \left\{ (x,y) \in \Omega :
f(x,y) - f \left( x_{0}, y_{0} \right) = A \left( x - x_{0}
\right) + B \left( y - y_{0} \right) \right\} .$$ It is not
difficult to see that since $\Omega $ is compact and $f$ is
$C^{1}$, and so continuous, it follows that $f(x,y) - f \left(
x_{0}, y_{0} \right) = A \left( x - x_{0} \right) + B \left( y -
y_{0} \right) $ is a closed condition and consequently $C(A,B)$ is
a closed subset of $\Omega $, which implies in its turn that
$C(A,B) \in K( \Omega )$. (See, for example, (3.15.1) and (3.17.3)
of [6].)
\\ \\
{\bf 1.3. Theorem.} If $$\left( A_{n}, B_{n} \right) \rightarrow
\left( \frac{ \partial f}{\partial x} \left( x_{0}, y_{0} \right)
, \frac{ \partial f}{\partial y} \left( x_{0}, y_{0} \right)
\right) $$ in ${\bf R}^{2}$ as $n \rightarrow \infty $, then
$$\limsup\limits_{n \rightarrow \infty } C \left( A_{n}, B_{n}
\right) \subseteq C \left( \frac{ \partial f}{\partial x} \left(
x_{0}, y_{0} \right) , \frac{ \partial f}{\partial y} \left(
x_{0}, y_{0} \right) \right) $$ in $K( \Omega )$.
\\ \\
{\bf Proof.} If $\left( x_{0}^{*}, y_{0}^{*} \right) \in
\limsup\limits_{n \rightarrow \infty } C \left( A_{n}, B_{n}
\right) $, then there exists a subsequence $\left( \left(
A_{n_{k}} , B_{n_{k}} \right) \right) _{k \in {\bf N}}$ of the
sequence $\left( \left( A_{n} , B_{n} \right) \right) _{n \in {\bf
N}}$ such that for any $k \in {\bf N}$, there exists $\left(
x_{n_{k}} , y_{n_{k}} \right) \in C \left( A_{n_{k}} , B_{n_{k}}
\right) $ with the property that $$\lim\limits_{k \rightarrow
\infty } \left( x_{n_{k}} , y_{n_{k}} \right) = \left( x_{0}^{*},
y_{0}^{*} \right) ,$$ so $$f \left( x_{n_{k}} , y_{n_{k}} \right)
- f \left( x_{0}, y_{0} \right) = A_{n_{k}} \left( x_{n_{k}} -
x_{0} \right) + B_{n_{k}} \left( y_{n_{k}} - y_{0} \right) $$ for
every $k \in {\bf N}$ and by passing to the limit as $k
\rightarrow \infty $ it follows that $$f \left( x_{0}^{*} ,
y_{0}^{*} \right) - f \left( x_{0}, y_{0} \right) = A \left(
x_{0}^{*} - x_{0} \right) + B \left( y_{0}^{*} - y_{0} \right) ,$$
so $\left( x_{0}^{*} , y_{0}^{*} \right) \in C \left( \frac{
\partial f}{\partial x} \left( x_{0}, y_{0} \right) , \frac{
\partial f}{\partial y} \left( x_{0}, y_{0} \right) \right) $.
\hfill $\bigtriangleup $
\\ \\
{\bf 1.4. Remark.} If ${\phi }(x) = \left\{ \begin{array}{ll} e^{ -
\frac{1}{x} } & \mbox{if $x>0$}
\\
0 & \mbox{if $x \leq 0$}
\end{array}
\right.$ and $\Omega = [-1,1]^{2}$, while $f(x,y) = {\phi }(x)$,
whenever $(x,y) \in \Omega $, then, since $\phi $ is $C^{1}$, it
follows that $f$ is $C^{1}$ and the tangent plane at $(0,0)$ is
the one with equation $z=0$, since $f(0,0)={\phi }(0)=0$, $\frac{
\partial f}{\partial x}(0,0)={\phi }'(0)=0$ and $\frac{
\partial f}{\partial y}(0,0)=0$. So
\begin{enumerate}
\item[ ]
$C(0,0) $
\item[ ]
$= \left\{ (x,y) \in [-1,1]^{2} : f(x,y) - f(0,0) = 0 \cdot (x-0)
+ 0 \cdot (y-0) \right\} $
\item[ ]
$= \left\{ (x,y) \in [-1,1]^{2} : {\phi }(x)=0 \right\} $
\item[ ]
$= [-1,0] \times [-1,1]$,
\end{enumerate}
while
\begin{enumerate}
\item[ ]
$C \left( \frac{1}{n} , 0 \right) $
\item[ ]
$= \left\{ (x,y) \in [-1,1]^{2} : f(x,y) - f(0,0) = \frac{1}{n}
\cdot (x-0) + 0 \cdot (y-0) \right\} $
\item[ ]
$= \left\{ (x,y) \in [-1,1]^{2} : {\phi }(x) = \frac{1}{n} x
\right\} $
\item[ ]
$= \left\{ (x,y) \in [-1,0) \times [-1,1] : 0 = \frac{1}{n} x
\right\} $
\item[ ]
$\cup \left\{ (x,y) \in \{ 0 \} \times [-1,1] : 0 = \frac{1}{n}
\cdot 0 \right\} $
\item[ ]
$\cup \left\{ (x,y) \in (0,1] \times [-1,1] : e^{ - \frac{1}{x} }
= \frac{1}{n} x \right\} $
\item[ ]
$= \{ 0 \} \times [-1,1]$
\item[ ]
$\cup \left\{ (x,y) \in (0,1] \times [-1,1] : e^{ - \frac{1}{x} }
= \frac{1}{n} x \right\} $
\end{enumerate}
and consequently $\limsup\limits_{n \rightarrow \infty } C \left(
\frac{1}{n} , 0 \right) $ is a proper subset of $C(0,0)$.
\\ \\
{\bf 1.5. Remark.} So, it follows that the tangent plane is the
limiting position of the normal vector of a secant plane.
\\ \\
{\bf 1.6. Definition.} Let $\Omega \in K \left( {\bf R}^{3} \right) $ be such that ${\Omega }^{ \circ } \neq \emptyset $ and let $f : \Omega \rightarrow {\bf R}$ be $C^{1}$, while $\left( x_{0}, y_{0}, z_{0} \right) \in {\Omega }^{ \circ }$. If $A$, $B$, $\Gamma $ are any real numbers and
\begin{enumerate}
\item[ ]
$P(A,B, \Gamma ) = \left\{ (x,y,z,t) \in {\bf R}^{4} : \right.$
\item[ ]
$\left. t - f \left( x_{0}, y_{0}, z_{0} \right) = A \left( x - x_{0} \right) + B \left( y - y_{0} \right) + \Gamma \left( z - z_{0} \right) \right\} $
\end{enumerate}
is any hyperplane in ${\bf R}^{4}$ that passes through $\left( x_{0}, y_{0}, z_{0}, f \left( x_{0}, y_{0}, z_{0} \right) \right) $ (see, for example, page 192 of [7] and page 295 of [8] and page 27 of [13]), then $P(A,B, \Gamma )$ intersects the hypersurface $Graph(f)$ whose equation is $$t = f(x,y,z)$$ exactly at the points $(x,y,z,t)$ in ${\bf R}^{4}$ that satisfy the equation $$f(x,y,z) - f \left( x_{0}, y_{0}, z_{0} \right) = A \left( x - x_{0} \right) + B \left( y - y_{0} \right) + \Gamma \left( z - z_{0} \right) $$ and let
\begin{enumerate}
\item[ ]
$H(A,B, \Gamma ) = \left\{ (x,y,z) \in \Omega : \right.$
\item[ ]
$\left. f(x,y,z) - f \left( x_{0}, y_{0}, z_{0} \right) = A \left( x - x_{0} \right) + B \left( y - y_{0} \right) + \Gamma \left( z - z_{0} \right) \right\} .$
\end{enumerate}
It is not difficult to see that since $\Omega $ is compact and $f$ is $C^{1}$, and so continuous, it follows that $$f(x,y,z) - f \left( x_{0}, y_{0}, z_{0} \right) = A \left( x - x_{0} \right) + B \left( y - y_{0} \right) + \Gamma \left( z - z_{0} \right) $$ is a closed condition and consequently $H(A,B, \Gamma )$ is a closed subset of $\Omega $, which implies in its turn that $H(A,B, \Gamma ) \in K( \Omega )$. (See, for example, (3.15.1) and (3.17.3) of [6].)
\\ \\
{\bf 1.7. Theorem.} If $$\left( A_{n}, B_{n}, {\Gamma }_{n} \right) \rightarrow \left( \frac{ \partial f}{\partial x} \left( x_{0}, y_{0}, z_{0} \right) , \frac{ \partial f}{\partial y} \left( x_{0}, y_{0}, z_{0} \right) , \frac{ \partial f}{\partial z} \left( x_{0}, y_{0}, z_{0} \right) \right) $$ in ${\bf R}^{3}$ as $n \rightarrow \infty $, then $$\limsup\limits_{n \rightarrow \infty } H \left( A_{n}, B_{n}, {\Gamma }_{n} \right) \subseteq H \left( \frac{ \partial f}{\partial x} \left( x_{0}, y_{0}, z_{0} \right) , \frac{ \partial f}{\partial y} \left( x_{0}, y_{0}, z_{0} \right) , \frac{ \partial f}{\partial z} \left( x_{0}, y_{0}, z_{0} \right) \right) $$ in $K( \Omega )$.
\\ \\
{\bf Proof.} If $\left( x_{0}^{*}, y_{0}^{*}, z_{0}^{*} \right) \in \limsup\limits_{n \rightarrow \infty } H \left( A_{n}, B_{n}, {\Gamma }_{n} \right) $, then there exists a subsequence $\left( \left( A_{n_{k}} , B_{n_{k}} , {\Gamma }_{n_{k}} \right) \right) _{k \in {\bf N}}$ of the sequence $\left( \left( A_{n} , B_{n} , {\Gamma }_{n} \right) \right) _{n \in {\bf N}}$ such that for any $k \in {\bf N}$, there exists $\left( x_{n_{k}} , y_{n_{k}} , z_{n_{k}} \right) \in H \left( A_{n_{k}} , B_{n_{k}} , {\Gamma }_{n_{k}} \right) $ with the property that $$\lim\limits_{k \rightarrow \infty } \left( x_{n_{k}} , y_{n_{k}} , z_{n_{k}} \right) = \left( x_{0}^{*}, y_{0}^{*}, z_{0}^{*} \right) ,$$ so
\begin{enumerate}
\item[ ]
$f \left( x_{n_{k}} , y_{n_{k}} , z_{n_{k}} \right) - f \left( x_{0}, y_{0}, z_{0} \right) $
\item[ ]
$= A_{n_{k}} \left( x_{n_{k}} - x_{0} \right) + B_{n_{k}} \left( y_{n_{k}} - y_{0} \right) + {\Gamma }_{n_{k}} \left( z_{n_{k}} - z_{0} \right) $
\end{enumerate}
for every $k \in {\bf N}$ and by passing to the limit as $k \rightarrow \infty $ it follows that
\begin{enumerate}
\item[ ]
$f \left( x_{0}^{*} , y_{0}^{*} , z_{0}^{*} \right) - f \left( x_{0}, y_{0}, z_{0} \right) $
\item[ ]
$= A \left( x_{0}^{*} - x_{0} \right) + B \left( y_{0}^{*} - y_{0} \right) + \Gamma \left( z_{0}^{*} - z_{0} \right) ,$
\end{enumerate}
so $\left( x_{0}^{*} , y_{0}^{*} , z_{0}^{*} \right) \in H \left( \frac{ \partial f}{\partial x} \left( x_{0}, y_{0}, z_{0} \right) , \frac{ \partial f}{\partial y} \left( x_{0}, y_{0}, z_{0} \right) , \frac{ \partial f}{\partial z} \left( x_{0}, y_{0}, z_{0} \right) \right) $. \hfill $\bigtriangleup $
\\ \\
{\bf 1.8. Remark.} If ${\phi }(x) = \left\{ \begin{array}{ll} e^{ - \frac{1}{x} } & \mbox{if $x>0$} \\ 0 & \mbox{if $x \leq 0$} \end{array} \right.$ and $\Omega = [-1,1]^{3}$, while $f(x,y,z) = {\phi }(x)$, whenever $(x,y,z) \in \Omega $, then, since $\phi $ is $C^{1}$, it follows that $f$ is $C^{1}$ and the tangent hyperplane at $(0,0,0)$ is the one with equation $t=0$, since $f(0,0,0)={\phi }(0)=0$, $\frac{ \partial f}{\partial x}(0,0,0)={\phi }'(0)=0$, $\frac{ \partial f}{\partial y}(0,0,0)=0$ and $\frac{ \partial f}{\partial z}(0,0,0)=0$. So
\begin{enumerate}
\item[ ]
$H(0,0,0) = \left\{ (x,y,z) \in [-1,1]^{3} : \right.$
\item[ ]
$\left. f(x,y,z) - f(0,0,0) = 0 \cdot (x-0)
+ 0 \cdot (y-0) + 0 \cdot (z-0) \right\} $
\item[ ]
$= \left\{ (x,y,z) \in [-1,1]^{3} : {\phi }(x)=0 \right\} $
\item[ ]
$= [-1,0] \times [-1,1]^{2}$,
\end{enumerate}
while
\begin{enumerate}
\item[ ]
$H \left( \frac{1}{n} , 0, 0 \right) = \left\{ (x,y,z) \in [-1,1]^{3} : \right.$
\item[ ]
$\left. f(x,y,z) - f(0,0,0) = \frac{1}{n} \cdot (x-0) + 0 \cdot (y-0) + 0 \cdot (z-0) \right\} $
\item[ ]
$= \left\{ (x,y,z) \in [-1,1]^{3} : {\phi }(x) = \frac{1}{n} x \right\} $
\item[ ]
$= \left\{ (x,y,z) \in [-1,0) \times [-1,1]^{2} : 0 = \frac{1}{n} x
\right\} $
\item[ ]
$\cup \left\{ (x,y,z) \in \{ 0 \} \times [-1,1]^{2} : 0 = \frac{1}{n} \cdot 0 \right\} $
\item[ ]
$\cup \left\{ (x,y,z) \in (0,1] \times [-1,1]^{2} : e^{ - \frac{1}{x} } = \frac{1}{n} x \right\} $
\item[ ]
$= \{ 0 \} \times [-1,1]^{2}$
\item[ ]
$\cup \left\{ (x,y,z) \in (0,1] \times [-1,1]^{2} : e^{ - \frac{1}{x} } = \frac{1}{n} x \right\} $
\end{enumerate}
and consequently $\limsup\limits_{n \rightarrow \infty } H \left( \frac{1}{n} , 0 , 0 \right) $ is a proper subset of $H(0,0,0)$.
\\ \\
{\bf 1.9. Remark.} So, it follows that the tangent hyperplane is the
limiting position of the normal vector of a secant hyperplane.

\section{One theorem in mechanics}

Our purpose in this section is to prove that every $C^{2}$ space
curve can be the solution of the principle of stationary action.
See Section 13 on page 59 of [1], where the principle is erroneously called of least action.
To this end, it is enough to prove the following:
\\ \\
{\bf 2.1. Theorem.} If $$L \left( t,x,y,z,x',y',z' \right) =
P(t,x,y,z) + Q_{1}(t,x)x' + Q_{2}(t,y)y' + Q_{3}(t,z)z',$$
whenever $\left( t,x,y,z,x',y',z' \right) \in {\bf R}^{7}$, where
\[
\left\{
\begin{array}{lllll}
\frac{ \partial P}{\partial x} = \frac{ \partial Q_{1}}{\partial
t}
\\ \\
\frac{\partial P}{\partial y} = \frac{ \partial Q_{2}}{\partial t}
\\ \\
\frac{ \partial P}{\partial z} = \frac{ \partial Q_{3}}{\partial
t}
\end{array}
\right.
\]
then given $- \infty <a<b< \infty $ and given any $C^{2}$ function
$${\bf r} : (a,b) \ni t \mapsto {\bf r}(t) = \left( x(t),y(t),z(t)
\right) \in {\bf R}^{3},$$ the function ${\bf r}$ is a solution of
the principle of stationary action with Lagrangian $L$, i.e.,
${\bf r}$ is a stationary point of the action integral $$I( {\bf
r} ) = \int_{a}^{b} L \left( t,{\bf r}(t),{\bf r}'(t) \right)
dt.$$

\noindent {\bf Proof.} By virtue of D.4 on page 152 of [15], it
follows that ${\bf r}$ must solve the system of equations
\[
\left\{
\begin{array}{lllll}
\frac{\partial L}{\partial x} - \frac{d}{dt} \left( \frac{\partial
L}{\partial x'} \right) = 0
\\ \\
\frac{\partial L}{\partial y} - \frac{d}{dt} \left( \frac{\partial
L}{\partial y'} \right) = 0
\\ \\
\frac{\partial L}{\partial z} - \frac{d}{dt} \left( \frac{\partial
L}{\partial z'} \right) = 0
\end{array}
\right.
\]
So, since $$\frac{\partial L}{\partial x} = \frac{\partial
P}{\partial x} + \frac{\partial Q_{1}}{\partial x} x'$$ and
$$\frac{\partial L}{\partial x'} = Q_{1}(t,x) \Rightarrow
\frac{d}{dt} \left( \frac{\partial L}{\partial x'} \right) =
\frac{\partial Q_{1}}{\partial t} + \frac{\partial Q_{1}}{\partial
x} x',$$ due to $$\frac{ \partial P}{\partial x} = \frac{ \partial
Q_{1}}{\partial t},$$ it follows that $$\frac{\partial L}{\partial
x} - \frac{d}{dt} \left( \frac{\partial L}{\partial x'} \right) =
0.$$ In addition, since $$\frac{\partial L}{\partial y} =
\frac{\partial P}{\partial y} + \frac{\partial Q_{2}}{\partial y}
y'$$ and $$\frac{\partial L}{\partial y'} = Q_{2}(t,y) \Rightarrow
\frac{d}{dt} \left( \frac{\partial L}{\partial y'} \right) =
\frac{\partial Q_{2}}{\partial t} + \frac{\partial Q_{2}}{\partial
y} y',$$ due to $$\frac{ \partial P}{\partial y} = \frac{ \partial
Q_{2}}{\partial t},$$ it follows that $$\frac{\partial L}{\partial
y} - \frac{d}{dt} \left( \frac{\partial L}{\partial y'} \right) =
0.$$ And once more, since $$\frac{\partial L}{\partial z} =
\frac{\partial P}{\partial z} + \frac{\partial Q_{3}}{\partial z}
z'$$ and $$\frac{\partial L}{\partial z'} = Q_{3}(t,z) \Rightarrow
\frac{d}{dt} \left( \frac{\partial L}{\partial z'} \right) =
\frac{\partial Q_{3}}{\partial t} + \frac{\partial Q_{3}}{\partial
z} z',$$ due to $$\frac{ \partial P}{\partial z} = \frac{ \partial
Q_{3}}{\partial t},$$ it follows that $$\frac{\partial L}{\partial
z} - \frac{d}{dt} \left( \frac{\partial L}{\partial z'} \right) =
0.$$ \hfill $\bigtriangleup $
\\ \\
{\bf 2.2. Theorem.} Keeping the notation as in {\bf 2.1.}, if
$$P(t,x,y,z) = P_{1}(t,x) + P_{2}(t,y) + P_{3}(t,z),$$ whenever
$(t,x,y,z) \in {\bf R}^{4}$, then the system of equations
\[
\left\{
\begin{array}{lllll}
\frac{ \partial P}{\partial x} = \frac{ \partial Q_{1}}{\partial
t}
\\ \\
\frac{\partial P}{\partial y} = \frac{ \partial Q_{2}}{\partial
t}
\\ \\
\frac{ \partial P}{\partial z} = \frac{ \partial Q_{3}}{\partial
t}
\end{array}
\right.
\]
is equivalent to the system of equations
\[
\left\{
\begin{array}{lllll}
\frac{ \partial P_{1}}{\partial x} = \frac{ \partial
Q_{1}}{\partial t}
\\ \\
\frac{\partial P_{2}}{\partial y} = \frac{ \partial
Q_{2}}{\partial t}
\\ \\
\frac{ \partial P_{3}}{\partial z} = \frac{ \partial
Q_{3}}{\partial t}
\end{array}
\right.
\]

\noindent {\bf Proof.} It is enough to notice that $$\frac{
\partial P}{\partial x} = \frac{ \partial P_{1}}{\partial x}$$ and
$$\frac{\partial P}{\partial y} = \frac{\partial P_{2}}{\partial
y}$$ and $$\frac{ \partial P}{\partial z} = \frac{ \partial
P_{3}}{\partial z}.$$ \hfill $\bigtriangleup $
\\ \\
{\bf 2.3. Theorem.} Keeping the notation as in {\bf 2.2.}, one can
solve easily
\[
\left\{
\begin{array}{lllll}
\frac{ \partial P_{1}}{\partial x} = \frac{ \partial
Q_{1}}{\partial t}
\\ \\
\frac{\partial P_{2}}{\partial y} = \frac{ \partial
Q_{2}}{\partial t}
\\ \\
\frac{ \partial P_{3}}{\partial z} = \frac{ \partial
Q_{3}}{\partial t}
\end{array}
\right.
\]

\noindent {\bf Proof.} By virtue of 4.1.5 on page 147 of [5],
since $$\frac{ \partial P_{1}}{\partial x} = \frac{
\partial Q_{1}}{\partial t},$$ the differential equation
$$P_{1}(t,x)dt + Q_{1}(t,x)dx = 0$$ has as solution the equation
$$u_{1}(t,x) = C,$$ where $C \in {\bf R}$, while $$u_{1}(t,x) =
\int_{t_{1}}^{t}P_{1}(t,x)dt + \int_{x_{1}}^{x}Q_{1} \left( t_{1},
x \right) dx,$$ so it is enough for $P_{1}$, $Q_{1}$ to be
$C^{1}$. In addition, by virtue of 4.1.5 on page 147 of [5], since
$$\frac{ \partial P_{2}}{\partial y} = \frac{
\partial Q_{2}}{\partial t},$$ the differential equation
$$P_{2}(t,y)dt + Q_{2}(t,y)dy = 0$$ has as solution the equation
$$u_{2}(t,y) = C,$$ where $C \in {\bf R}$, while $$u_{2}(t,y) =
\int_{t_{2}}^{t}P_{2}(t,y)dt + \int_{y_{2}}^{y}Q_{2} \left( t_{2},
y \right) dy,$$ so it is enough for $P_{2}$, $Q_{2}$ to be
$C^{1}$. And once more, by virtue of 4.1.5 on page 147 of [5],
since $$\frac{ \partial P_{3}}{\partial z} = \frac{
\partial Q_{3}}{\partial t},$$ the differential equation
$$P_{3}(t,z)dt + Q_{3}(t,z)dz = 0$$ has as solution the equation
$$u_{3}(t,z) = C,$$ where $C \in {\bf R}$, while $$u_{3}(t,z) =
\int_{t_{3}}^{t}P_{3}(t,z)dt + \int_{z_{3}}^{z}Q_{3} \left( t_{3},
z \right) dz,$$ so it is enough for $P_{3}$, $Q_{3}$ to be
$C^{1}$. \hfill $\bigtriangleup $
\\ \\
{\bf 2.4. Remark.} In order to do infinitely many examples of $L$'s, do the following : Take three $C^{2}$ functions $f(t,x)$, $g(t,y)$, $h(t,z)$ and set $$L = \frac{ \partial f }{ \partial t } + \frac{ \partial g }{ \partial t } + \frac{ \partial h }{ \partial t } + \frac{ \partial f }{ \partial x } x' + \frac{ \partial g }{ \partial y } y' + \frac{ \partial h }{ \partial z } z'.$$ For example, $$L = \left( 2t + 3x^{2}t \right) + \left( 2ty \right) + \left( -2z \right) + \left( 3t^{2}x \right) x' + \left( t^{2} - y^{2} \right) y' + \left( z - 2t \right) z'.$$ Since for any such $L$, the equations of motion in Lagrangian mechanics give as solution every $C^{2}$ space
curve and since $L=T-U$, where $T$ is kinetic energy and $U$ is potential energy, when kinetic energy $T$ is given, it follows that potentials $U=T-L$ with $L$ as above do not exist. See Section 13 on page 59 of [1].

\end{document}